\documentclass{amsart}
\usepackage{hyperref}

\usepackage{amsrefs,amsmath,amsfonts,amssymb,graphicx,latexsym,color}

\renewcommand{\dim}{n}

\newcommand{\threed}{{\mathbb R}^\dim}

\newcommand{\R}{{\mathbb R}}
\newcommand{\sph}{{\mathbb S}^{\dim-1}}

\newcommand{\eqdef}{\overset{\mbox{\tiny{def}}}{=}}

\newcommand{\weirdQ}{\rho}

\newcommand{\pZ}{p^{0}}
\newcommand{\qZ}{q^{0}}
\newcommand{\pZprime}{p^{\prime0}}
\newcommand{\qZprime}{q^{\prime0}}

\newcommand{\pM}{p^{\mu}}
\newcommand{\qM}{q^{\mu}}
\newcommand{\pMprime}{p^{\prime\mu}}
\newcommand{\qMprime}{q^{\prime\mu}}

\newcommand{\angleN}{k}

\newcommand{\relMOM}{\varrho}

\newtheorem{theorem}{Theorem}
\newtheorem{corollary}[theorem]{Corollary}

\theoremstyle{definition}

\newtheorem*{rem}{Remark}
\newtheorem{definition}[theorem]{Definition}

\begin{document}
\title[Coordinates in the relativistic Boltzmann theory]{Coordinates in the relativistic Boltzmann theory}

\author[R. M. Strain]{Robert M. Strain}
\address{
University of Pennsylvania, Department of Mathematics, David
Rittenhouse Lab, 209 South 33rd Street, Philadelphia, PA 19104}
\email{strain at math.upenn.edu}
\urladdr{http://www.math.upenn.edu/~strain/}
\thanks{R.M.S. was partially supported by the NSF grant DMS-0901463.}
\keywords{Relativity, Boltzmann, 
 collisional Kinetic Theory, Kinetic Theory}
\subjclass[2000]{Primary: 76P05; Secondary: 83A05}

\setcounter{tocdepth}{1}

\begin{abstract}
It is often the case in mathematical analysis that solving an open problem can be facilitated by finding a new set of coordinates  which may illumniate the known difficulties. In this article, we illustrate how to derive an assortment coordinates in which to represent the relativistic Boltzmann collision operator.  We show the equivalence between some known representations \cite{MR1211782,MR635279}, and others which seem to be new.  One of these representations has been used recently to solve several open problems in
 \cite{strainNEWT,strainSOFT,gsRVMB,ssHilbert}.
\end{abstract}

\maketitle
\tableofcontents

 \medskip
 \dedicatory{Dedicated to the memory of 
Naoufel ben Abdallah
and
Carlo Cercignani}

\section{Introduction and main result}

The relativistic Boltzmann equation can be written as 
\begin{equation*}
\pM \partial_\mu F = \mathcal{C}(F,F).
\end{equation*}
In this expression the collision operator \cite{MR1026740,MR635279} is given by 
\begin{equation}
\mathcal{C}(f,h) = \int_{\threed} \frac{dq}{\qZ}
\int_{\threed}\frac{dq^\prime}{\qZprime}
\int_{\threed}\frac{dp^\prime}{\pZprime}
W(p, q | p^\prime, q^\prime) [f(p^{\prime})h(q^{\prime})-f(p)h(q)].  
\notag
\end{equation}
The transition rate, $W(p, q | p^\prime, q^\prime)$, from \cite{MR1026740} in $\dim$-dimensions ($n\ge 2$) is denoted 
\begin{equation}
W(p, q | p^\prime, q^\prime) = \frac{1}{2} 
 \left(\frac{\relMOM}{2}\right)^{3-\dim} ~ s ~ \sigma(\relMOM, \theta) ~
\delta^{(1+\dim)}(\pM+\qM-\pMprime-\qMprime),  \label{transition}
\end{equation}
where $\sigma(\relMOM, \theta)$ is the differential cross-section which is a measure of the interactions between particles.  
This is an important model for fast moving
particles.  Standard
references in relativistic Kinetic theory include \cite{MR1898707,MR635279,MR1379589,MR0088362,stewart}. 
In this paper we give a complete reduction of the collision integrals for the operator $\mathcal{C}(f,h)$, deriving several sets of coordinates for the particles momentum, some old and some new.
The rest of our notation is given after the following historical discussion.

\subsection{Historical discussion}
The relativistic Boltzmann equation is the primary model in relativistic collisional Kinetic theory.  In the next few paragraphs, we will provide a short review of the mathematical theory of this equation.  We mention a few books on relativistic Kinetic theory as for instance \cite{MR1898707,MR635279,MR1379589,MR0088362,stewart}.  

Carlo Cercignani and collaborators wrote several research works on relativity and the Boltzmann equation including
\cite{MR700073,MR1770447,MR767047,MR1706733,MR1745845,MR1810262,MR1857127,MR1898707}; 
 many of these works were in collaboration with 
Gilberto Medeiros Kremer \cite{MR1706733,MR1745845,MR1810262,MR1857127,MR1898707}.  
Together \cite{MR1706733} they showed  that the original proof of the formula for the ``summational invariants'', due to Boltzmann, also applies in the special relativistic case.  In addition Cercignani-Kremer wrote a book on the relativistic Boltzmann equation \cite{MR1898707}.
This follows the book by  Cercignani-Illner-Pulvirenti \cite{MR1307620} on the Newtonian Boltzmann equation.  The author learned about the physics of the relativistic Boltzmann equation from the books of
Cercignani-Kremer \cite{MR1898707}
and
de Groot-van Leeuwen-van Weert \cite{MR635279}.

Now early results include those of Lichnerowicz and Marrot \cite{MR0004796}, who wrote down the full relativistic Boltzmann equation, including collisional effects, in 1940.  Later, in 1967, Bichteler \cite{MR0213137} showed that the relativistic Boltzmann equation has a local solution.  Then Dudy{\'n}ski and Ekiel-Je{\.z}ewska \cite{MR933458}, in 1988, proved that the linearized equation admits unique solutions in $L^2$.  Afterwards, Dudy{\'n}ski \cite{MR1031410} studied the long time and small-mean-free-path limits of these solutions.

In the context of large data global in time weak solutions, the theory of DiPerna-Lions \cite{MR1014927} renormalized solutions was extended  to the relativistic Boltzmann equation in 1992 by Dudy{\'n}ski and Ekiel-Je{\.z}ewska \cite{MR1151987}.  This result uses the causality of the relativistic Boltzmann equation 
\cite{MR818441,MR841735}.  
Results on the regularity of the gain term are given in \cite{MR1402446,MR1480243};
 the strong $L^1$ compactness is studied by  Andr{\'e}asson \cite{MR1402446}.
These are generalizations of Lions \cite{MR1284432} result in the non-relativistic case. 
Further developments on renormalized weak solutions can be found in 
\cite{MR1714446,MR1676150}.

Let us point out studies of the Newtonian limit \cite{MR2098116,strainNEWT} for the Boltzmann equation.  
We mention theories of unique global in time solutions with initial data that is ``near Vacuum'' as in \cite{MR2217287,strainNEWT,MR2459827}.  Andr{\'e}asson, Calogero and Illner \cite{MR2102321} prove that there can be blow-up
in the presence of only the  gain term.
Notice the study of the collision map
 and the pre-post collisional change of variables from \cite{MR1105532}.  
Then \cite{MR2543323} provides uniform $L^2$-stability estimates for the relativistic
   Boltzmann equation.  
 Now there is a mathematically rigorous result connecting the relativistic Euler equations to the relativistic Boltzmann equation via the Hilbert expansion as in \cite{ssHilbert}.

Next we discuss several results for initial data that is nearby the relativistic Maxwellian.
In 1993, Glassey and Strauss \cite{MR1211782} first proved the global existence and uniqueness of smooth solutions  on the torus $\mathbb{T}^3$.  They also established exponential convergence to Maxwellian. Their assumptions on the differential cross-section, $\sigma(\relMOM,\theta)$, covered the case of hard potentials.  
 In 1995, they extended that result to 
the whole space case \cite{MR1321370} where the convergence rate is polynomial.  Subsequent results with reduced restrictions 
on the cross-section were  proven in \cite{MR2249574}, using the energy method from for instance \cite{MR2000470}; these results also apply to the hard potentials.  Recently, the author \cite{strainSOFT} has developed a weighted 
$L^\infty(\mathbb{T}^3_x\times \mathbb{R}^3_p)$ theory of global in time unique solutions for the soft potentials, and rapid convergence to Maxwellian, under the general physical assumption proposed in \cite{DEnotMSI,MR933458}.

We point out results on global existence of unique smooth solutions which are initially close to the relativistic Maxwellian for the 
relativistic Landau-Maxwell system \cite{MR2100057}, and then for the relativistic Landau \cite{MR2289548} equation as well.  Further \cite{MR2514726} proves the smoothing effects for relativistic Landau-Maxwell system.  
And \cite{MR2593052} proves time decay rates in the whole space for the relativistic Boltzmann equation (with hard potentials) and the relativistic Landau equation as well.

\subsection{Notation}
In this section we define several notations which will be used throughout the article.
A relativistic particle has momentum $p=(p^1, \ldots, p^\dim)\in\threed$, with 
its energy defined by $\pZ=\sqrt{c^2+|p|^2}$ where $|p|^2 \eqdef p\cdot p$.  
Here $c$ denotes the speed of light.  
We use the standard Euclidean dot product: $p\cdot q \eqdef \sum_{i=1}^{\dim} p^i q^i$.  
As is customary we write $\pM = (\pZ, p)$ where $\pM$ also denotes the $\mu$-th element of $(\pZ, p)$.
In general, Latin (spatial) indices $i,j,$ etc., take values in $\{1,\ldots,\dim\},$ while Greek indices $\kappa, \lambda, \mu, \nu,$ etc., take on the values $\{0,1,\ldots,\dim\}.$
Indices are raised and lowered as usual with the Minkowski metric $g_{\mu \nu}$ and its inverse $g^{\mu \nu}$,
where
$(g_{\mu\nu}) \eqdef \text{diag}(-1 ~ 1 ~ \cdots ~ 1)$ is an $(1+\dim) \times (1+\dim)$ matrix. 
In other words
$p_\mu = g_{\mu\nu} p^\nu$.  
We furthermore use the Einstein convention of implicit summation over repeated indices.
The Lorentz inner product  is then given by
$$
\pM q_\mu = \pM g_{\mu\nu} q^\nu \eqdef -\pZ \qZ+p\cdot q.
$$
Notice 
$
\pM p_\mu = -c^2.
$
Now 
$\partial_\mu = (c^{-1}\partial_t, \nabla_x)$ so the streaming term of the relativistic Boltzmann equation is 
$
\pM \partial_\mu 
=
c^{-1}p_0\partial_t + p\cdot \nabla_x,
$
where $\nabla_x$ is the spatial gradient.

Conservation of momentum and energy for elastic collisions is expressed as
\begin{equation}
\pM+\qM=\pMprime+\qMprime.
\label{collisionalCONSERVATION}
\end{equation}
These conservation laws are enforced by the $1+n$ delta functions in \eqref{transition}.
The angle $\theta$ in the Boltzmann collision operator \eqref{transition} is then defined by
\begin{equation}
\cos\theta
\eqdef
(\pM - \qM) (p_\mu^\prime -q_\mu^\prime)/\relMOM^2.
\label{angle}
\end{equation}
Note that this angle is well defined under \eqref{collisionalCONSERVATION}, see  \cite[Lemma 3.15.3]{MR1379589}.
Here the relative momentum, $\relMOM = \relMOM(\pM, \qM)$,  is denoted
\begin{gather}
\relMOM 
\eqdef
\sqrt{(\pM-\qM) (p_\mu-q_\mu)}
=
\sqrt{2(-\pM q_\mu-c^2)}
=
\sqrt{2(p^0 q^0 - p\cdot q-c^2)}.
\label{gDEFINITION}
\end{gather}
Furthermore  the quantity $s= s(\pM, \qM)$ is defined as
\begin{eqnarray}
s 
\eqdef
-(\pM+\qM) (p_\mu+q_\mu)
=
2\left( -\pM q_\mu+c^2\right)\ge 0.
\label{sDEFINITION}
\end{eqnarray}
This is in other words
$
s = 2(p^0 q^0 - p\cdot q+c^2).
$
Notice that $s=\relMOM^2+4c^2$.   To proceed, we will quickly review the Lorentz transformations.

\subsection*{A digression on Lorentz transformations}
In this subsection, we discuss a few elementary aspects of Lorentz transformations which will be useful throughout the rest of this paper.  
Let $\Lambda$ be a $(1+\dim)\times (1+\dim)$ matrix (of real numbers) denoted by
$$
\Lambda=(\Lambda_{\ \nu}^{\mu})_{0\le \mu,\nu\le \dim}.
$$
For the basics of Lorentz transformations, we refer to \cite{MR1898707} and \cite{weinbergBK}.  

\begin{definition}\label{rcop:LTdef}
$\Lambda$ is a (proper) Lorentz transformation if $\mbox{det}(\Lambda) = 1$ and 
$$
\Lambda_{\ \mu}^{\kappa} ~ g_{\kappa \lambda} ~ \Lambda_{\ \nu}^{\lambda} = g_{\mu \nu}, \qquad (\mu, \nu = 0,1,\ldots,\dim).
$$
This implies the following invariance:
$
p^\kappa q_\kappa = p^\kappa g_{\kappa\lambda } q^{\lambda} 
= 
(\Lambda^{\kappa}_{~\mu}  p^{\mu})g_{\kappa\lambda } (\Lambda^{\lambda}_{~\nu}  q^{\nu} ).
$
\end{definition}

We will use the notation $\Lambda$ to exclusively denote a Lorentz transformation.
Any such $\Lambda$ is invertible and
the inverse matrix is denoted $\Lambda^{-1} = (\Lambda_{\mu}^{~\nu})_{0\le \mu,\nu\le \dim}$
so that 
$
\Lambda_{~\kappa}^{\nu}\Lambda_{\mu}^{~\kappa}= \delta_{~\mu}^{\nu},
$
where 
$
\delta_{~\mu}^{\nu}
$
is the standard Kronecker delta which is unity when the indices are equal and zero otherwise.
It follows from Definition \ref{rcop:LTdef} that
\begin{equation}
\Lambda_{\mu}^{~\nu}= g^{\nu \lambda }~  \Lambda_{~ \lambda }^{\kappa} ~  g_{\kappa \mu}.
\label{Linverse}
\end{equation}
Definition \ref{rcop:LTdef} further implies that  $(\Lambda^{-1})^{\mu}_{~\nu}  = \Lambda^{~\mu}_{\nu}$ 
is a Lorentz transformation.


In this paper, we are exclusively concerned with (proper)  Lorentz transformations $\Lambda$, depending on $\pM$ and $\qM$, which further satisfy
\begin{equation}
\Lambda^{\mu}_{~\nu}(p^{\nu}+q^{\nu})=(\sqrt{s}, 0, \ldots, 0).
\label{lw:lorentzS}
\end{equation}
We recall that $s$ is defined in \eqref{sDEFINITION}.  
Note that there are several Lorentz transformations which will satisfy \eqref{lw:lorentzS}.  In particular, in dimension $\dim =3$, three such Lorentz transformations were given in \cite[Appendix A]{strainNEWT}.

The condition \eqref{lw:lorentzS} means that you are mapping the particle momentum to the ``center of momentum'' frame $p+q=0$ (this is sometimes called alternatively the center of mass frame).  Notice that \eqref{lw:lorentzS} further implies 
\begin{equation}
\Lambda^{0}_{~\mu}(\pM-\qM)=0.
\label{lorentzSnull}
\end{equation}
This follows quickly from the fact that \eqref{lw:lorentzS} implies 
$
\Lambda^{i}_{~\mu}\pM= -\Lambda^{i}_{~\mu}\qM
$
for any $i\in \{1, \ldots, \dim\}$.  
Then Definition \ref{rcop:LTdef} allows us to further deduce \eqref{lorentzSnull}
as a result of the mass shell condition: 
$
-c^2 = p^\kappa p_\kappa
= 
(\Lambda^{\kappa}_{~\mu}  p^{\mu})g_{\kappa\lambda } (\Lambda^{\lambda}_{~\nu}  p^{\nu} ).
$

We give a standard example of these mappings.

\subsubsection*{The Boost matrix} 
The most common Lorentz transformation is probably the Boost matrix.  
Given 
$
v=(v^1,\ldots,v^\dim)\in\threed,
$  
the $(1+\dim)\times (1+\dim)$  Boost matrix is
$$
\Lambda_{b}
\eqdef
\left(
\begin{array}{cc}
\weirdQ & -\weirdQ v
\\
-\weirdQ v & {\bf 1}_n+(\weirdQ-1)\frac{v \otimes v}{|v|^2} 
\end{array}
\right),
$$
where $\weirdQ=(1-|v|^2)^{-1/2}$ and ${\bf 1}_n$ is the $\dim \times \dim$ identity matrix.  Notice that $\Lambda_{b}$ has only $\dim$ free parameters.  
Our goal is to choose $v$ such that $\Lambda_{b}$ satisfies (\ref{lw:lorentzS}).  Let
$$
v\eqdef \frac{p+q}{\pZ+\qZ}, \quad \weirdQ\eqdef \frac{\pZ+\qZ}{\sqrt{s}}.
$$
Plugging these choices into $\Lambda_{b}$ above, we obtain that
\begin{equation}
\Lambda_{b}
=
\left(
\begin{array}{cc}
\frac{\pZ+\qZ}{\sqrt{s}} & -\frac{p+q}{\sqrt{s}}
\\
-\frac{p+q}{\sqrt{s}} & {\bf 1}_n+(\weirdQ-1)\frac{(p+q) \otimes (p+q)}{|p+q|^2} 
\end{array}
\right).
\label{boostLT}
\end{equation}
By a direct calculation this example satisfies (\ref{lw:lorentzS}).

\subsection{Main results}
We will now state our main results in the language of the notation just introduced.  
For this we use the operator $\mathcal{Q}(f,h) \eqdef \mathcal{C}(f,h)\frac{c}{\pZ}$ as
\begin{equation}
\mathcal{Q}(f,h)
=
\frac{c}{\pZ}\int_{\threed} \frac{dq}{\qZ}
\int_{\threed}\frac{dq^\prime}{\qZprime}
\int_{\threed}\frac{dp^\prime}{\pZprime}
W(p, q | p^\prime, q^\prime) [f(p^{\prime})h(q^{\prime})-f(p)h(q)].  
\label{collisionTQ}
\end{equation}
The main point is to carry out the reduction of the number of integrations in this expression by evaluating the $1+\dim$ delta functions from \eqref{transition}.  In the literature, there are two approaches to performing this goal \cite{MR1211782,MR635279}.  One \cite{MR1211782} uses algebraic manipulation of polynomials; this results in a representation 
(Theorem \ref{GStransformationPROP}) that is well known in the mathematics literature and has been widely used (at least in dimension $\dim = 3$).  
The other representation \cite{MR635279} (from Corollary \ref{corBOOST}  and Theorem \ref{CMtransformationPROP}), which starts by using the change of variables as in \eqref{lw:lorentzS}, is sketched in physics texts 
but seems to be hard to locate in the mathematics literature on the relativistic Boltzmann equation; this approach was developed (as in Corollary \ref{corBOOST} below in dimension $\dim = 3$) in the authors thesis \cite{strainPHD}.  Additionally the representations of the collision operator \eqref{collisionTQ} from Corollary \ref{corBOOST} have only recently been used to solve several open problems on the relativistic Boltzmann equation 
\cite{strainNEWT,strainSOFT,gsRVMB,ssHilbert}.

For this reason, we are motivated at this time to write down a complete mathematical proof of these different representations.  
Our results improve upon those given previously \cite{MR1211782,MR635279} in the following ways.  We prove these representations in $\dim$ dimensions with $\dim \ge 2$.  Note that the representation in Theorem \ref{GStransformationPROP} is dimensionally dependent; there is an additional term \eqref{kernelGS} which is not present when $\dim =3$.   Furthermore it is discussed in the physics paper \cite{MR1026740} that there are physical situations in which the Boltzmann equation may be of interest in dimensions other than $\dim =3$.
We also prove the exact formula for the post-collisional energies as in \eqref{pcE}.
We give the precise expressions for the angles in \eqref{angleDEFtheta}, 
\eqref{angleDEF},
and
\eqref{angleGS}
which do not seem to have been previously computed.  What the author finds most interesting is that we can show that there are several alternative expressions for the post-collisional momentum and energy, as in 
\eqref{rcop:postCOLLvel} and \eqref{rcop:postCOLLvelENER}, one for each Lorentz transformation satisfying \eqref{lw:lorentzS}.
We hope that the availability of these additional alternative representations may be useful to future investigations in the relativistic Boltzmann theory; in particular we observe in Corollary \ref{comboCOR} the equivalence of these different representations.  This equivalence has been crucial to our recent proof (joint with Yan Guo) of the global in time stability for the relativistic Vlasov-Maxwell-Boltzmann system \cite{gsRVMB} with near Maxwellian initial conditions.  

To begin we state the center of momentum \eqref{lw:lorentzS} reduction.

\begin{theorem}
\label{CMtransformationPROP} 
(Center of momentum reduction). 
Recall \eqref{collisionTQ} and %
\eqref{transition}. For any suitable integrable function $G: \threed\times%
\threed\times\threed\times\threed\to \R$, it holds that
\begin{multline}
\frac{c}{\pZ \qZ} \int_{\threed} \frac{dq^\prime}{\qZprime}\int_{%
\threed}\frac{dp^\prime}{\pZprime} W(p, q | p^\prime, q^\prime) G(p,
q, p^{\prime}, q^{\prime}) 
\\
= \int_{\sph} d\omega ~ v_{\o } ~ \sigma(\relMOM, \theta) ~ 
G(p, q, p^{\prime}, q^{\prime}),  
\label{rcop:gainS}
\end{multline}
where $\omega = (\omega^1, \ldots, \omega^\dim) \in \sph$ and
$v_{\o }=v_{\o }(p,q)$ is the M\o ller velocity given by 
\begin{equation}
v_{\o}= v_{\o}(p,q) \eqdef 
\frac{c}{4}
\frac{ \relMOM\sqrt{s}}{\pZ \qZ}.  \label{moller}
\end{equation}
The angle $\theta$ in this expression is defined by \eqref{angleDEFtheta} and \eqref{angleDEF}.  
The post-collisional momentum and energy above are defined by \eqref{rcop:postCOLLvel} and \eqref{rcop:postCOLLvelENER} respectively.  
\end{theorem}


The angle in the reduced expression in Theorem \ref{CMtransformationPROP} is defined by
\begin{equation}
\cos\theta = \frac{\angleN}{|\angleN|}\cdot \omega,
\label{angleDEFtheta}
\end{equation}
where $\angleN\in\threed$ is given as
\begin{equation}
\angleN^i \eqdef \Lambda^{i}_{~\mu} (\pM-\qM)
=
\Lambda^{i}_{~0} (\pZ-\qZ)
+
\Lambda^{i}_{~j} (p^j-q^j),
\quad (i=1,\ldots, \dim),
\label{angleDEF}
\end{equation}
Moreover the post-collisional momentum satisfy ($i=1,\ldots,\dim$)
\begin{equation}
\begin{split}
p^{\prime i}
        &=
L^{i0}  \frac{\sqrt{s}}{2}
- 
L^{i j} \omega_j \frac{\relMOM}{2},
\\
q^{\prime i}
    &=
L^{i0}  \frac{\sqrt{s}}{2}
+
L^{i j} \omega_j \frac{\relMOM}{2}.
\end{split}
\label{rcop:postCOLLvel}
\end{equation}
Here we use the tensor $L^{\mu\kappa} \eqdef -g^{\mu \lambda }~  \Lambda_{~ \lambda }^{\kappa}$.
Furthermore, the energies are
\begin{equation}
\begin{split}
\pZprime
        &=
L^{00}  \frac{\sqrt{s}}{2}
- 
L^{0 j} \omega_j \frac{\relMOM}{2},
\\
\qZprime
&=
L^{00}  \frac{\sqrt{s}}{2}
+ 
L^{0 j} \omega_j \frac{\relMOM}{2}.
\end{split}
\label{rcop:postCOLLvelENER}
\end{equation}
As usual, in each of these expressions, we implicitly sum over $j \in\{ 1,\ldots, \dim \}$.  Note that
the formulae above hold for any Lorentz transformation $\Lambda$ satisfying only \eqref{lw:lorentzS}.

Now \eqref{transition}, \eqref{collisionTQ} and \eqref{rcop:gainS} together imply that
$$
\mathcal{Q}(f,h)
=
 \int_{\threed} dq
\int_{\sph} d\omega ~ v_{\o } ~ \sigma(\relMOM, \theta) ~  [f(p^{\prime})h(q^{\prime})-f(p)h(q)],
$$
where the relevant quantities are defined above as in \eqref{angleDEFtheta}, \eqref{angleDEF}, \eqref{rcop:postCOLLvel} and \eqref{rcop:postCOLLvelENER}.  

 If we use the Lorentz boost \eqref{boostLT} satisfying \eqref{lw:lorentzS} we have the following simplification.

\begin{corollary}\label{corBOOST}
In the particular case of the Boost matrix \eqref{boostLT}, the post-collisional momentum \eqref{rcop:postCOLLvel}
 in the expressions above can be written precisely as
\begin{equation}
\begin{split}
p^\prime&=\frac{p+q}{2}+\frac{\relMOM}{2}\left(\omega+(\weirdQ-1)(p+q)\frac{%
(p+q)\cdot \omega}{|p+q|^2}\right), \\
q^\prime&=\frac{p+q}{2}-\frac{\relMOM}{2}\left(\omega+(\weirdQ-1)(p+q)\frac{%
(p+q)\cdot \omega}{|p+q|^2}\right),  \label{postCOLLvelCMsec2}
\end{split}%
\end{equation}
where $\weirdQ \eqdef (\pZ+\qZ)/\sqrt{s}$.  Furthermore, the post-collisional 
energies are given by
\begin{equation}
\begin{split}
\pZprime
        &=
   \frac{\pZ+\qZ}{2}
+  
\frac{\relMOM}{2\sqrt{s}}
~(p+q)\cdot   \omega,
\\
\qZprime
&=
   \frac{\pZ+\qZ}{2}
-
\frac{\relMOM}{2\sqrt{s}}
~(p+q)\cdot   \omega.
\end{split}
\notag
\end{equation}
Additionally in the angle \eqref{angleDEFtheta}, the vector \eqref{angleDEF} can be simplified to
\begin{equation}
\angleN
=
-\frac{p+q}{\sqrt{s}} (\pZ-\qZ)
+
(p-q)
+
(\rho -1)(p+q)\frac{(p+q)\cdot (p-q)}{|p+q|^2}.
\notag
\end{equation}
In these formula we use $s$ from \eqref{sDEFINITION}.
 \end{corollary}

Now we turn to the expression given by Glassey-Strauss in \cite{MR1211782}.

\begin{theorem}
\label{GStransformationPROP} 
(Glassey-Strauss reduction). 
Alternatively 
\begin{multline*}
\frac{c}{\pZ \qZ} \int_{\threed} \frac{dq^\prime}{\qZprime}\int_{%
\threed}\frac{dp^\prime}{\pZprime} W(p, q | p^\prime, q^\prime) G(p,
q, p^{\prime}, q^{\prime})  \notag \\
= \int_{\sph}d\omega ~ \frac{s\sigma(\relMOM,\theta)}{\pZ\qZ} ~ B_\dim(p,q,\omega) ~
G(p, q, p^{\prime}, q^{\prime}),
\end{multline*}
where $B_\dim(p,q,\omega)$ is given by \eqref{kernelGS} and $(p^\prime, q^\prime)$
on the right are given by \eqref{postCOLLvelGS}.  The angle $\theta$ is also defined by \eqref{angleGS}.
\end{theorem}

In the above reduction, we consider the expression
\begin{equation}
B(p,q,\omega )\eqdef c\frac{%
(p^0+q^0)^{2}p^0q^0\left\vert \omega \cdot \left( \frac{p}{p^0}-%
\frac{q}{q^0}\right) \right\vert }{\left[ (p^0+q^0)^{2}-(\omega \cdot
\lbrack p+q])^{2}\right] ^{2}}.  \label{kernelGSthr}
\end{equation}%
Then the kernel in Theorem \ref{GStransformationPROP} is given by
\begin{equation}
B_\dim (p,q,\omega )\eqdef 
B(p,q,\omega )  
\left(
A(p,q,\omega )  
\right)^{\dim-3},
\label{kernelGS}
\end{equation}%
where
$$
A (p,q,\omega ) 
\eqdef
\frac{(p^0+q^0)p^0q^0
\left\vert \omega \cdot \left( \frac{p}{p^0}-\frac{q}{q^0}\right) \right\vert }{ \left[ (p^0+q^0)^{2}-(\omega \cdot
\lbrack p+q])^{2}\right]} \frac{4}{\relMOM}.
$$
Notice that $B_3(p,q,\omega ) = B(p,q,\omega )$ when the dimension is $\dim = 3$.

For the reduction in Theorem \ref{GStransformationPROP}, the post-collisional momentum are given by
\begin{equation}
\begin{split}
p^{\prime } &=p+a(p,q,\omega )\omega,
 \\
q^{\prime } &=q-a(p,q,\omega )\omega,
\end{split}
\label{postCOLLvelGS}
\end{equation}
where
\begin{equation*}
a(p,q,\omega )=\frac{2(\pZ+\qZ)\{\omega \cdot (\pZ q-\qZ p)\} }{(p^0+q^0)^{2}-(\omega \cdot \lbrack p+q])^{2}}.
\end{equation*}
And the energies can be expressed 
as
$\pZprime=\pZ+N^0$ and
$\qZprime=\qZ-N^0$
with 
\begin{equation}
N^0\eqdef \frac{2 \omega\cdot(p+q)\{\omega \cdot (\pZ q-\qZ p)\}}{(p^0+q^0)^{2}-(\omega \cdot \lbrack p+q])^{2}}.
\label{pcE}
\end{equation}
These formula clearly satisfy the collisional conservations 
\eqref{collisionalCONSERVATION}.  
The angle \eqref{angle} in Theorem \ref{GStransformationPROP}, can then be reduced to
\begin{equation}
\cos\theta
=
1
-
\frac{8}{\relMOM^2}
\frac{
\{\omega \cdot (\pZ q-\qZ p)\}^2}{ (p^0+q^0)^{2}-(\omega \cdot \lbrack p+q])^{2} }.
\label{angleGS}
\end{equation}
Evidently, in this case $\theta$ is not simply a dot product of a unit vector with $\omega$.

Moreover, assuming the collisions are elastic as in \eqref{collisionalCONSERVATION}, we have the invariance:
$$
 \omega \cdot \left( \pZ  q- \qZ  p\right)  
=
 \omega \cdot \left( \qZprime p^{\prime}- \pZprime q^{\prime}\right).
$$
Therefore, for fixed $\omega\in \sph$, 
$
B_\dim(p,q,\omega)=B_\dim(p^\prime ,q^\prime,\omega).
$
Then this kernel is invariant under  pre-post change of variables $(p', q') \to (p,q)$ as in \cite{MR1105532}.

In this case \eqref{transition}, \eqref{collisionTQ} and Theorem \ref{GStransformationPROP} together illustrate that
$$
\mathcal{Q}(f,h)
=
 \int_{\threed} dq
\int_{\sph} d\omega ~ \frac{s\sigma(\relMOM,\theta)}{\pZ\qZ} ~ B_\dim(p,q,\omega) ~  [f(p^{\prime})h(q^{\prime})-f(p)h(q)],
$$
where the quantities inside the integral are defined as in \eqref{angleGS}, \eqref{postCOLLvelGS}, \eqref{kernelGS} and \eqref{kernelGSthr}.

We mention now the following useful corollary.

\begin{corollary}\label{comboCOR}
Combining Theorem \ref{CMtransformationPROP} and Theorem \ref{GStransformationPROP} yields
\begin{gather*}
\int_{\sph} d\omega ~ v_{\o } ~ \sigma(\relMOM, \theta) ~ G(p, q, p^\prime ,q^\prime) =
\int_{\sph}d\omega ~ \frac{s\sigma(\relMOM,\theta)}{\pZ\qZ} ~ B_\dim(p,q,\omega) ~ 
G(p,q, p^\prime ,q^\prime),
\end{gather*}
where the angle and post-collisional momentum on the left are defined as in Theorem \ref{CMtransformationPROP}, and the similar expressions on the right are given as in Theorem \ref{GStransformationPROP}.
 \end{corollary}

Now in the variables \eqref{postCOLLvelGS}, it was computed in \cite[Theorem 1]{MR1105532} that the mapping $(p,q) \to (p', q')$ has the following Jacobian 
\begin{equation*}
\left| \frac{\partial (p^{\prime },q^{\prime })}{\partial (p,q)} \right| = \frac{p_{0}^{\prime }q_{0}^{\prime }}{p_{0}q_{0}}.
\end{equation*}
Note also that when $(p,q) \to (p', q')$ then additionally $(p', q') \to (p,q)$.  We conclude
\begin{multline*}
\int_{\threed} dp \int_{\threed} dq
\int_{\sph}d\omega ~ \frac{s\sigma(\relMOM,\theta)}{\pZ\qZ} ~ B_\dim(p,q,\omega) ~ 
G(p,q, p^\prime ,q^\prime)
\\
=
\int_{\threed} dp \int_{\threed} dq
\int_{\sph}d\omega ~ \frac{s\sigma(\relMOM,\theta)}{\pZ\qZ} ~ B_\dim(p,q,\omega) ~ 
G(p^\prime ,q^\prime, p,q).
\end{multline*}
This holds with the variables \eqref{postCOLLvelGS} which are used in Theorem \ref{GStransformationPROP}.  However iterating this formula and using Corollary \ref{comboCOR} we additionally observe that
\begin{multline}
\int_{\threed} dp \int_{\threed} dq
\int_{\sph}d\omega
 ~ v_{\o } ~ \sigma(\relMOM, \theta) ~ G(p, q, p^\prime ,q^\prime)
\\
=
\int_{\threed} dp \int_{\threed} dq
\int_{\sph}d\omega 
 ~ v_{\o } ~ \sigma(\relMOM, \theta) ~ 
G(p^\prime ,q^\prime, p,q).
\end{multline}
The variables in these integrals are those from \eqref{rcop:postCOLLvel} which are used in Theorem \ref{CMtransformationPROP}.

This article is organized as follows.  In the next Section \ref{secAPP:reduction} we will prove Theorem \ref{CMtransformationPROP}, from which we conclude Corollary \ref{corBOOST}.  Then in Section \ref{secAPP:reduction2} we will prove Theorem \ref{GStransformationPROP}.  We remark that Corollary \ref{comboCOR} will be an important part of \cite{gsRVMB}.

\section{Center of Momentum reduction of the Collision Integrals}\label{secAPP:reduction}

In this section we will prove Theorem \ref{CMtransformationPROP} and in particular \eqref{rcop:gainS}.  
To this end we consider the following integral
\begin{multline}
I \eqdef
\int_{\threed} \frac{dq^\prime}{\qZprime}\int_{\threed}\frac{dp^\prime}{\pZprime}
~ s \sigma(\relMOM, \theta)  ~ \delta^{(1+\dim)}(\pM+\qM-\pMprime-\qMprime)
G(\pM,\qM,\pMprime,\qMprime)
\\
=
\int_{\threed} \frac{dq^\prime}{\qZprime}\int_{\threed}\frac{dp^\prime}{\pZprime}
~ s \sigma(\relMOM, \theta)  ~ \delta^{(1+\dim)}(\pM+\qM-\pMprime-\qMprime)
G(p, q, p^{\prime}, q^{\prime}),
\label{IintR}
\end{multline}
where $\theta$ is defined by \eqref{angle}.  Above and in the following, for convenience, we write 
$G(\pM,\qM,\pMprime,\qMprime)=G(p, q, p^{\prime}, q^{\prime})$ when there is no opportunity for confusion. 

Our goal will be to prove that 
\begin{equation}
I 
=2^{2-\dim}\int_{\sph} d\omega ~ \relMOM^{n-2} \sqrt{s} ~ \sigma(\relMOM, \theta) ~ 
G(p, q, p^{\prime}, q^{\prime}).  
\label{IintRgoal}
\end{equation}
Above we use the definitions \eqref{angleDEFtheta}, \eqref{angleDEF}, \eqref{rcop:postCOLLvel}, and \eqref{rcop:postCOLLvelENER}.  Notice then, by \eqref{transition} and \eqref{moller}, 
that this implies \eqref{rcop:gainS}.

\subsection*{Integral reduction}   We now reduce the collision integrals in the center of momentum system.
We first claim that $I$ from \eqref{IintR} can be written as
$$
I  =
\int_{\threed} \frac{dq^\prime}{\qZprime}\int_{\threed}\frac{dp^\prime}{\pZprime}
s \sigma(\relMOM, \theta) \delta^{(1+\dim)}(\Lambda^{\mu}_{~\nu}(p^{\nu}+q^{\nu})-\pMprime-\qMprime)
G(\pM,\qM,\tilde{p}^{\prime\mu},\tilde{q}^{\prime\mu}).
$$
Here $\Lambda^{\mu}_{~\nu}$ should satisfy \eqref{lw:lorentzS}.
In the expression above from \eqref{Linverse} we have
\begin{equation}
\tilde{p}^{\prime\mu}
\eqdef
\Lambda^{~\mu}_{\nu} p^{\prime \nu}
=
g^{\mu \lambda }~  \Lambda_{~ \lambda }^{\kappa} ~  g_{\kappa \nu} ~ p^{\prime \nu},
\quad 
\text{and}
\quad 
\tilde{q}^{\prime\mu}
= 
g^{\mu \lambda }~  \Lambda_{~ \lambda }^{\kappa} ~  g_{\kappa \nu} ~ q^{\prime \nu}.
\label{tildePRIMEv}
\end{equation}
This holds because $\frac{dq^\prime}{\qZprime}$ and $\frac{dp^\prime}{\pZprime}$ are Lorentz invariant measures and 
$$
\delta^{(1+\dim)}(\Lambda^{\mu}_{~\nu}(p^{\nu}+q^{\nu}-p^{\prime\nu}-q^{\prime\nu}) )
=
\delta^{(1+\dim)}(\pM+\qM-\pMprime-\qMprime ).
$$
But notice that the claim is not true unless the angle $\theta$ from \eqref{angle} is redefined as
$$
\cos\theta
=
\{\Lambda^{\nu}_{~\mu} (\pM - \qM) \}(p_\nu^\prime -q_\nu^\prime)/\relMOM^2
=\frac{\angleN}{|\angleN|}\cdot \frac{(p^\prime-q^\prime)}{\relMOM},
$$
where $k$ is defined in \eqref{angleDEF}.  We have also employed the following calculation
$$
\relMOM = 
\sqrt{(\pM-\qM) (p_\mu-q_\mu)}
=
\sqrt{\{\Lambda^{\mu}_{~\kappa} (p^\kappa-q^\kappa)\}g_{\mu\nu}\{ \Lambda^{\nu}_{~\lambda} (p^\lambda-q^\lambda)\}}
=|\angleN|.
$$
For that we used \eqref{lorentzSnull}.  Note also, as a result of \eqref{lw:lorentzS}, we can deduce (from the delta function in $I$) that $p^\prime+q^\prime=0$ which further yields $\pZprime = \qZprime$.

Then the integration over $q^\prime$ can be carried out immediately and we obtain
$$
I  =
\frac 12 \int_{\threed} \frac{dp^\prime}{(\pZprime)^2}
s \sigma(\relMOM, \theta) \delta\left(\frac 12 \sqrt{s}-\pZprime\right)
G(\pM,\qM,\tilde{p}^{\prime\mu},\tilde{q}^{\prime\mu}).
$$
Here now, with $p^\prime=-q^\prime$, the angle is
$$
\cos\theta= 2\frac{\angleN}{|\angleN|}\cdot \frac{p^\prime}{\relMOM}.
$$
Next change to polar coordinates as $p^\prime=|p^\prime|\omega$ with $\omega\in \sph$ 
and
$$
dp^\prime=|p^\prime|^{\dim-1} d|p^\prime| d\omega.
$$
We use the following calculation to compute the delta function
\begin{eqnarray*}
\frac12 \sqrt{s}-\pZprime
&=&
\frac12 \sqrt{s}-\sqrt{c^2+|p^\prime|^2}
=
\frac{\frac14 s-(c^2+|p^\prime|^2)}{\frac 12 \sqrt{s}+\sqrt{c^2+|p^\prime|^2}}
\\
&=&
\frac{\frac14  \relMOM^2-|p^\prime|^2}{\frac 12 \sqrt{s}+\sqrt{c^2+|p^\prime|^2}}.
\end{eqnarray*}
We have just employed  $s= \relMOM^2 + 4c^2$ from \eqref{gDEFINITION} and \eqref{sDEFINITION}.  Thus 
$$
\frac12 \sqrt{s}-\pZprime
=
\frac{(\frac12 \relMOM -|p^\prime|)(\frac12 \relMOM+|p^\prime|)}{\frac 12 \sqrt{s}+\sqrt{c^2+|p^\prime|^2}}.
$$
Now if $\frac12 \relMOM=|p^\prime|$ then $s= \relMOM^2 + 4c^2$ grants that
$
\sqrt{c^2+|p^\prime|^2}=\frac 12 \sqrt{s}.
$
Therefore,
$$
 \delta\left(\frac12 \sqrt{s}-\pZprime\right)
 =
\frac{\frac 12 \sqrt{s}+\sqrt{c^2+|p^\prime|^2}}{ \frac 12 \relMOM+|p^\prime|}
\delta\left(\frac12 \relMOM-|p^\prime|\right)
 =
\frac{\sqrt{s}}{\relMOM}\delta\left(\frac12 \relMOM-|p^\prime|\right).
$$
We plug this in to see that
\begin{multline}
\notag
I =
\frac 12 \int_0^\infty d|p^\prime| \int_{\sph}d\omega\frac{|p^\prime|^{\dim-1} }{(\pZprime)^2}
\frac{s^{3/2}}{\relMOM} \sigma(\relMOM, \theta) \delta\left(\frac12 \relMOM-|p^\prime|\right)
G(\pM,\qM,\tilde{p}^{\prime\mu},\tilde{q}^{\prime\mu})
\\
=
\left(\frac{1}{2}\right)^{\dim }  \int_{\sph}d\omega ~ \frac{\relMOM^{\dim-1} }{(\frac{1}{2}\sqrt{s})^2}
\frac{s^{3/2}}{\relMOM} \sigma(\relMOM, \theta) 
G(\pM,\qM,p^{\prime\mu},q^{\prime\mu})
\\
=
\left(\frac{1}{2}\right)^{\dim - 2}  \int_{\sph}d\omega ~ \relMOM^{\dim-2} 
s^{1/2} \sigma(\relMOM, \theta) 
G(\pM,\qM,p^{\prime\mu},q^{\prime\mu}).
\end{multline}
Now consider $\bar{\omega}^\mu \eqdef (\pZprime, |p^\prime|\omega)=\frac 12(\sqrt{s}, \relMOM\omega)$
and 
$
\tilde{\omega}^\mu
\eqdef
(\pZprime, -|p^\prime|\omega)=\frac 12(\sqrt{s}, -\relMOM\omega).
$
Then in the integral for $I$ above, we replace \eqref{tildePRIMEv} with
\begin{equation}
\pMprime
=
g^{\mu \lambda }~  \Lambda_{~ \lambda }^{\kappa} ~  g_{\kappa \nu} ~ \bar{\omega}^\nu,
\quad
\qMprime
=
g^{\mu \lambda }~  \Lambda_{~ \lambda }^{\kappa} ~  g_{\kappa \nu} ~ \tilde{\omega}^\nu.
\label{postCcom}
\end{equation}
We calculate below that these are \eqref{rcop:postCOLLvel} and \eqref{rcop:postCOLLvelENER}.
Furthermore the angle is given as in \eqref{angleDEFtheta}.  
We thus obtain \eqref{IintRgoal}, (\ref{rcop:gainS}) and Theorem \ref{CMtransformationPROP}.

\subsection*{Post-collisional momentum}
We now calculate the post-collisional momentum and their associated energies \eqref{rcop:postCOLLvel} and \eqref{rcop:postCOLLvelENER}.
We use \eqref{postCcom}
to compute 
$$
\pMprime
=
-  L^{\mu\kappa}  g_{\kappa \nu}\bar{\omega}^\nu
=
-  L^{\mu\kappa}  \bar{\omega}_\kappa
=
L^{\mu0}  \frac{\sqrt{s}}{2}
-  
L^{\mu j} \omega_j \frac{\relMOM}{2},
$$
where  $L^{\mu\kappa} \eqdef -g^{\mu \lambda }~  \Lambda_{~ \lambda }^{\kappa}$.
Similarly
$$
\qMprime
=
-  L^{\mu\kappa}  \tilde{\omega}_\kappa
=
L^{\mu0}  \frac{\sqrt{s}}{2}
+ 
L^{\mu j} \omega_j \frac{\relMOM}{2}.
$$
These establish that \eqref{postCcom} is equal to \eqref{rcop:postCOLLvel} and \eqref{rcop:postCOLLvelENER}.

Notice that, using  \eqref{postCcom}, we recover the identities of conservation for elastic collisions (\ref{collisionalCONSERVATION}) from \eqref{Linverse} and \eqref{lw:lorentzS}  as follows
$$
\Lambda^{\mu}_{~\nu}(p^{\nu}+q^{\nu})
=
(\sqrt{s}, 0, \ldots, 0)
=  
\bar{\omega}^\mu
+
\tilde{\omega}^\mu
=
\Lambda^{\mu}_{~\nu}(p^{\prime\nu}+q^{\prime\nu}).
$$
The last equality follows from \eqref{postCcom} since
$$
\Lambda^{\mu}_{~\nu}(p^{\prime\nu}+q^{\prime\nu})
=  
\Lambda^{\mu}_{~\nu} 
~ g^{\nu \lambda }~  \Lambda_{~ \lambda }^{\kappa} ~  g_{\kappa \alpha} ~\left( \bar{\omega}^\alpha +  \tilde{\omega}^\alpha\right)
=
\bar{\omega}^\mu
+
\tilde{\omega}^\mu,
$$
which itself is implied by
$
\Lambda^{\mu}_{~\nu} 
g^{\nu \lambda }~  \Lambda_{~ \lambda }^{\kappa}
=
g^{\mu \kappa};
$
the previous statement is a consequence of Definition \ref{rcop:LTdef}.  To obtain \eqref{collisionalCONSERVATION} we apply the inverse, as in \eqref{Linverse}, to observe
$$
\pMprime + \qMprime 
=  
\Lambda^{~\mu}_{\kappa}
\Lambda^{\kappa}_{~\nu}(p^{\prime\nu}+q^{\prime\nu})
=
\Lambda^{~\mu}_{\kappa}
\Lambda^{\kappa}_{~\nu}(p^{\nu}+q^{\nu})
=
p^{\mu}+q^{\mu}.
$$
This completes the proof of Theorem \ref{CMtransformationPROP}.

\section{Glassey-Strauss reduction of the Collision Integrals}\label{secAPP:reduction2}

The goal of this section is to prove Theorem \ref{GStransformationPROP}.  This reduction was given by Glassey \& Strauss \cite{MR1211782} in dimension $\dim = 3$ and without the presence of the arguments, $(\relMOM, \theta)$, in the differential cross section, $\sigma(\relMOM, \theta)$.  The reduction below is essentially similar to \cite{MR1211782} in that we perform the examination of the roots of polynomials; it is also however slightly different from \cite{MR1211782} in that we integrate the radial variable, $r$, below over the entire real line $\R$ (rather than $r\ge 0$).  We also observe a dimensionally dependent factor $A(p,q,\omega)$ as in \eqref{kernelGS} which is not present when $\dim = 3$.  We are further able to compute the formula for the angle $\theta$ in \eqref{angleGS} and establish the equivalence of the representations as in Corollary \ref{comboCOR}.

We will use the expression for $I$ in \eqref{IintR}.  We claim that
\begin{equation}
I=  \frac{2}{c} \left(\frac{\relMOM}{2}\right)^{\dim-3}  \int_{\sph} d\omega 
~ s\sigma(\relMOM,\theta) ~
B(p,q,\omega )  
\left(
A(p,q,\omega )  
\right)^{\dim-3}
G(p, q, p^\prime, q^\prime),
\label{claimGS}
\end{equation}
where $(p^\prime, q^\prime)$ are defined above as in \eqref{postCOLLvelGS}.
In this case, using \eqref{IintR}, \eqref{transition}, \eqref{kernelGSthr} and \eqref{kernelGS}, 
we obtain Theorem \ref{GStransformationPROP} subject to \eqref{claimGS}.

To reduce the number of integrals in \eqref{IintR}, we split $I = \frac 12 I + \frac 12 I$.  
  Letting $q^\prime =p+q-p^\prime$ we can immediately remove $\dim$ of the delta functions in the first copy of $I$.  Next translate $p^\prime \to p+p^\prime$ so that $q^\prime \to q-p^\prime$.  Then switch to spherical coordinates as $p^\prime = r\omega$, $d p^\prime = r^{\dim-1} drd\omega$ where $r\in [0,\infty)$ and $\omega \in \sph$.  We obtain
$$
\frac 12I=  \frac 12\int_{0}^\infty \int_{\sph}  \frac{r^{\dim-1} dr d\omega }{\pZprime \qZprime} ~ s\sigma(\relMOM,\theta) ~
\delta(\pZprime+\qZprime-\pZ-\qZ)
G(p, q, p^\prime, q^\prime),
$$
where $p^\prime=p+r\omega$ and $q^\prime=q-r\omega$.  

Now consider the other $\frac 12I$.  Let $p^\prime =p+q-q^\prime$ to remove $\dim$ of the delta functions.  Next translate $q^\prime \to q+q^\prime$ so that $p^\prime \to p-q^\prime$.  And  switch to spherical coordinates as $q^\prime = r\omega$, 
$d q^\prime = r^{\dim-1} drd\omega$ where $r\in [0,\infty)$ and $\omega \in \sph$.  Then 
$$
\frac 12I=  \frac 12\int_{0}^\infty \int_{\sph} \frac{r^{\dim-1} dr d\omega }{\pZprime \qZprime}
~ s\sigma(\relMOM,\theta) ~
\delta(\pZprime+\qZprime-\pZ-\qZ)
G(p, q, p^\prime, q^\prime),
$$
where $p^\prime=p-r\omega$ and $q^\prime=q+r\omega$.  Further change variables $r\to -r$ so that 
$$
\frac 12I=  \frac 12\int_{-\infty}^0 \int_{\sph} \frac{|r|^{\dim-1} dr d\omega }{\pZprime \qZprime}
~ s\sigma(\relMOM,\theta) ~
\delta(\pZprime+\qZprime-\pZ-\qZ)
G(p, q, p^\prime, q^\prime),
$$
where now $p^\prime=p+r\omega$ and $q^\prime=q-r\omega$.

We combine the last two splittings to conclude that 
$$
I=\frac 12I+\frac 12I=  \frac 12\int_{-\infty}^{+\infty} \int_{\sph} \frac{|r|^{\dim-1} dr d\omega }{\pZprime \qZprime}
~ s\sigma(\relMOM,\theta) ~
\delta(\pZprime+\qZprime-\pZ-\qZ)
G(p, q, p^\prime, q^\prime),
$$
where $p^\prime=p+r\omega$ and $q^\prime=q-r\omega$.  Now the angle \eqref{angle} satisfies
\begin{multline}
\cos\theta
=
\frac{-(\pZ - \qZ)(\pZprime - \qZprime)+(p-q)\cdot (p' - q')}{\relMOM^2}
\\
=
\frac{-(\pZ - \qZ)(\pZprime - \qZprime)+|p-q|^2
+2r
(p-q)\cdot \omega}{\relMOM^2}.
\label{angle2}
\end{multline}
We will return to this expression below.

We now focus on the argument of the delta function.  For $\lambda_i>0$ ($i=1,2$), we use the identity 
$\delta(\lambda_1 - \lambda_2) = 2\lambda_1 \delta(\lambda_1^2 - \lambda_2^2)$ and \eqref{collisionalCONSERVATION} to see that
\begin{gather*}
\delta(\pZprime+\qZprime-\pZ-\qZ)
=
2(\pZ+\qZ)\delta((\pZprime+\qZprime)^2-(\pZ+\qZ)^2)
\\
=
2(\pZ+\qZ)~ \delta(2\pZprime \qZprime-\{(\pZ+\qZ)^2-(\pZprime)^2-(\qZprime)^2\}).
\end{gather*}
If $\{(\pZ+\qZ)^2-(\pZprime)^2-(\qZprime)^2\}>0$ then this is
\begin{gather*}
=
8(\pZ+\qZ)\pZprime \qZprime ~ \delta(4(\pZprime)^2 (\qZprime)^2-\{(\pZ+\qZ)^2-(\pZprime)^2-(\qZprime)^2\}^2).
\end{gather*}
If $\{(\pZ+\qZ)^2-(\pZprime)^2-(\qZprime)^2\}<0$, then the delta function is zero.
Now we write the argument of the last delta function above as
\begin{multline*}
p(r)=-(\pZ+\qZ)^4+2(\pZ+\qZ)^2\{(\pZprime)^2+(\qZprime)^2\}
\\
-\{(\pZprime)^2+(\qZprime)^2\}^2
+4(\pZprime)^2 (\qZprime)^2
\\
=-(\pZ+\qZ)^4+2(\pZ+\qZ)^2\{(\pZprime)^2+(\qZprime)^2\}
-\{(\pZprime)^2-(\qZprime)^2\}^2.
\end{multline*}
Plugging in $p^\prime=p+r\omega$ and  $q^\prime=q-r\omega $ we observe that 
$$
(\pZprime)^2-(\qZprime)^2
=
|p+r\omega|^2 -|q-r\omega|^2 = (\pZ)^2 -(\qZ)^2+2 r \omega \cdot (p+q).
$$
This means that $p(r)$ is quadratic in $r$.  Moreover,
\begin{equation*}
\begin{split}
p(0)&=-(\pZ+\qZ)^4+2(\pZ+\qZ)^2\{(\pZ)^2+(\qZ)^2\}
-\{(\pZ)^2-(\qZ)^2\}^2
\\
&=-(\pZ+\qZ)^4+2(\pZ+\qZ)^2\{(\pZ)^2+(\qZ)^2\}
-(\pZ+\qZ)^2(\pZ-\qZ)^2
\\
&=-(\pZ+\qZ)^4+(\pZ+\qZ)^2\{(\pZ)^2+(\qZ)^2\}
+2\pZ\qZ(\pZ+\qZ)^2=0.
\end{split}
\end{equation*}
 We conclude that $p(r)=4D_1 r^2-8D_2 r$ for some $D_1, D_2\in\mathbb{R}$.  
 
 We will now determine $D_1, D_2$.  Expanding
$$
(\pZprime)^2+(\qZprime)^2
=
2c^2 +|p+r\omega|^2 +|q-r\omega|^2 = (\pZ)^2+(\qZ)^2+2 r \omega \cdot (p-q)+2r^2.
$$
We plug the last few calculations into $p(r)$ to write it out in terms of $r$ as 
\begin{equation*}
\begin{split}
p(r)
&=-(\pZ+\qZ)^4+2(\pZ+\qZ)^2\{(\pZ)^2+(\qZ)^2+2 r \omega \cdot (p-q)+2r^2\}
\\
&\quad
-\{(\pZ)^2 -(\qZ)^2+2 r \omega \cdot (p+q)\}^2
\\
&=-(\pZ+\qZ)^4+2(\pZ+\qZ)^2\{(\pZ)^2+(\qZ)^2+2 r \omega \cdot (p-q)+2r^2\}
\\
&\quad
-\{(\pZ)^2 -(\qZ)^2\}^2-4 r^2\{ \omega \cdot (p+q)\}^2-4r\{ \omega \cdot (p+q)\}\{(\pZ)^2 -(\qZ)^2\}.
\end{split}
\end{equation*}
Rearrainging the terms
\begin{equation*}
\begin{split}
p(r)
&=
4\{(\pZ+\qZ)^2-\{\omega\cdot (p+q)\}^2\}r^2 
\\
&\quad
+
4\left\{(\pZ+\qZ)^2\{\omega \cdot (p-q)\}
-\{ \omega \cdot (p+q)\}\{(\pZ)^2 -(\qZ)^2\}\right\}r
\\
&\quad
-(\pZ+\qZ)^4+2(\pZ+\qZ)^2\{(\pZ)^2+(\qZ)^2\}
-\{(\pZ)^2-(\qZ)^2\}^2
\\
&=
4\{(\pZ+\qZ)^2-\{\omega\cdot (p+q)\}^2\}r^2 
\\
&\quad
+
4\left\{(\pZ+\qZ)^2\{\omega \cdot (p-q)\}
-\{ \omega \cdot (p+q)\}\{(\pZ)^2 -(\qZ)^2\}\right\}r.
\end{split}
\end{equation*}
Equivalently we have
\begin{equation*}
\begin{split}
D_1&=\{(\pZ+\qZ)^2-\{\omega\cdot (p+q)\}^2\}, 
\\
2D_2 &= -(\pZ+\qZ)^2\{\omega \cdot (p-q)\}
+\{ \omega \cdot (p+q)\}\{(\pZ)^2 -(\qZ)^2\}.
\end{split}
\end{equation*}
We further calculate $D_2$ as
\begin{equation*}
\begin{split}
2D_2 
&= -((\pZ)^2+(\qZ)^2+2\pZ\qZ)\{\omega \cdot (p-q)\}
+\{ \omega \cdot (p+q)\}\{(\pZ)^2 -(\qZ)^2\}
\\
&= 
w\cdot p\left\{-2(\qZ)^2-2\pZ\qZ \right\}
+\omega \cdot q \left\{2(\pZ)^2+2\pZ\qZ \right\}
\\
&=
2\{(\pZ+\qZ)\omega \cdot (\pZ q-\qZ p)\}
\\
&=
2(\pZ+\qZ)\pZ\qZ\left\{\omega \cdot \left(\frac{q}{\qZ}-\frac{p}{\pZ}\right)\right\}.
\end{split}
\end{equation*}
Thus, $p(r)=4D_1 r^2-8D_2 r$ with these definitions.

We plug this formulation for $p(r)$ into the full integral to obtain
$$
I=  \frac{1}{2}\int_{-\infty}^{+\infty} \int_{\sph} |r|^{\dim-1} dr d\omega 
~ s\sigma(\relMOM,\theta) ~
8(\pZ+\qZ)
\delta(4D_1 r^2-8D_2 r)
G(p, q, p^\prime, q^\prime).
$$
Then we extract the factor $4D_1$ from the delta function to obtain
$$
I=  \int_{-\infty}^{+\infty} \int_{\sph} |r|^{\dim-1} dr d\omega 
~ s\sigma(\relMOM,\theta) ~
\frac{(\pZ+\qZ)}{D_1}
\delta(r\{r-2D_2 /D_1\})
G(p, q, p^\prime, q^\prime).
$$
We use the identity 
$\delta(r\{r-2D_2 /D_1\})
=\left| \frac{D_1}{2D_2}\right|\{\delta(r)+\delta(r-2D_2 /D_1)\}$.  The first delta function drops out because of the $|r|^{\dim-1}$ factor ($\dim \ge 2$).  
We thus have
$$
I=  \int_{-\infty}^{+\infty} \int_{\sph} |r|^{\dim -1} dr d\omega 
~ s\sigma(\relMOM,\theta) ~
\frac{(\pZ+\qZ)}{2| D_2|}
\delta(r-2D_2 /D_1)
G(p, q, p^\prime, q^\prime).
$$
We have used $D_1\ge 2$ \cite{MR1105532}; note $D_2$ can be negative.  
Evaluate the delta function
$$
I=  2^{\dim -2}  \int_{\sph} d\omega 
~ s\sigma(\relMOM,\theta) ~
\left(\frac{|D_2|^{\dim - 2}}{D_1^{\dim - 1}}\right)
(\pZ+\qZ)
G(p, q, p^\prime, q^\prime).
$$
Notice that from \eqref{kernelGS} and \eqref{kernelGSthr} we have
$$
2^{\dim -2} \left(\frac{|D_2|^{\dim - 2}}{D_1^{\dim - 1}}\right)
(\pZ+\qZ)
=
\frac{2}{c}\left(\frac{\relMOM}{2}\right)^{\dim-3}
B(p,q,\omega )  
\left(
A(p,q,\omega )  
\right)^{\dim-3}. 
$$
This grants us \eqref{claimGS}, and thereby Theorem \ref{GStransformationPROP}.
Furthermore, the angle \eqref{angle2} becomes
\begin{multline}
\cos\theta
=
\frac{-(\pZ - \qZ)(\pZprime - \qZprime)+|p-q|^2
+2r
(p-q)\cdot \omega}{\relMOM^2}
\\
=
\frac{-(\pZ - \qZ)\frac{\{(\pZprime)^2 - (\qZprime)^2\}}{\pZ + \qZ}+|p-q|^2
+4\frac{D_2}{D_1}
(p-q)\cdot \omega}{\relMOM^2}
\\
=
\frac{-(\pZ - \qZ)\frac{\{(\pZ)^2 -(\qZ)^2+2 r \omega \cdot (p+q)\}}{\pZ + \qZ}+|p-q|^2
+4\frac{D_2}{D_1}
(p-q)\cdot \omega}{\relMOM^2}.
\notag
\end{multline}
Above we have used \eqref{collisionalCONSERVATION} and $r = 2D_2/D_1$. 
 Further calculations yield
\begin{multline}
\cos\theta
=
\frac{-(\pZ - \qZ)^2
-\frac{\pZ - \qZ}{\pZ + \qZ}~ 4\frac{D_2}{D_1}\omega \cdot (p+q)
+|p-q|^2
+4\frac{D_2}{D_1}
(p-q)\cdot \omega}{\relMOM^2}
\\
=
1
+
4\frac{D_2}{D_1}
\frac{
-\frac{\pZ - \qZ}{\pZ + \qZ}~ \omega \cdot (p+q)
+
(p-q)\cdot \omega}{\relMOM^2}
\\
=
1
+
4D_2
\frac{
-2\frac{\omega \cdot (\pZ q-\qZ p)}{(\pZ+\qZ)}}{D_1 ~ \relMOM^2}
\\
=
1
-
8
\frac{
\{\omega \cdot (\pZ q-\qZ p)\}^2}{D_1 ~ \relMOM^2}.
\label{angle3}
\end{multline}
Note \eqref{angle3} is exactly \eqref{angleGS}.

\subsection*{Post-collisional energy}  We show in this part how to compute the post-collisional energy \eqref{pcE}.  The point is to write $\pZprime=\pZ+N^0$ and $\qZprime=\qZ-N^0$, and also recall from 
$\pZprime=\sqrt{c^2 + |p^\prime|^2}$ and
$\qZprime=\sqrt{c^2 + |q^\prime|^2}$ that one can obtain
$$
-c^2 = \pMprime p_{\mu}^\prime
= -(\pZ+N^0)^2 + |p^\prime|^2
= -(\pZ+N^0)^2 + |p + a(p,q,\omega)\omega|^2.
$$
Similarly,
$$
-c^2 = \qMprime q_{\mu}^\prime
= -(\qZ-N^0)^2 + |q - a(p,q,\omega)\omega|^2.
$$
By expanding both of these expressions, and subtracting the result, one obtains
$$
2(\pZ + \qZ) N^0 - 2 \omega\cdot (p+q) a(p,q,\omega) = 0.
$$
Solving this equality for $N^0$ establishes \eqref{pcE}. 
A related computation (but for a different purpose) can be found in Cercignani-Kremer \cite[Section 1.4.3]{MR1898707}.
\hfill {\bf Q.E.D.}

\begin{rem}
We point out that there is unfortunately a misprint in our recent paper \cite{strainNEWT} if the dimension $\dim \ge 2$ is not $\dim = 3$.  In particular the transition rate $W$ at the top of the paper \cite{strainNEWT}, should be replaced by the transition rate from \eqref{transition} (and \cite{MR1026740}).  The main difference between the two is the factor $\relMOM^{3-\dim}$, which is unity when $\dim = 3$.  
Furthermore, the expression (1.9) in \cite{strainNEWT} for the collision operator is correct when $\dim = N = 3$ ($N$ is the notation for the dimension used in \cite{strainNEWT}), but otherwise the kernel 
$B(p,q,\omega)$ in \cite[(1.9)]{strainNEWT} needs to be replaced by $B_\dim(p,q,\omega)$ in \eqref{kernelGS}.
In other words the factor $\left(A(p,q,\omega)\right)^{\dim -3}$ from \eqref{kernelGS} is missing from \cite[(1.9)]{strainNEWT}.
We point out that this factor in the expression \cite[(1.9)]{strainNEWT} does not affect the main theorems of \cite{strainNEWT}.  
Furthermore the condition 
in \cite[(2.7)]{strainNEWT} on the collisional cross section
is written as  $0\le \gamma < -3$ (which is empty); this  \cite[(2.7)]{strainNEWT} 
should be $0\le \gamma < N$.  
\end{rem}



\begin{bibdiv}
\begin{biblist}


\bib{MR1402446}{article}{
    author={Andr{\'e}asson, H{\aa}kan},
     title={Regularity of the gain term and strong $L\sp 1$ convergence to
            equilibrium for the relativistic Boltzmann equation},
   journal={SIAM J. Math. Anal.},
    volume={27},
      date={1996},
    number={5},
     pages={1386\ndash 1405},
      issn={0036-1410},
    review={MR1402446 (97e:76077)},
}

\bib{MR2102321}{article}{
    author={Andr{\'e}asson, H{\aa}kan},
    author={Calogero, Simone},
    author={Illner, Reinhard},
     title={On blowup for gain-term-only classical and relativistic
            Boltzmann equations},
   journal={Math. Methods Appl. Sci.},
    volume={27},
      date={2004},
    number={18},
     pages={2231\ndash 2240},
      issn={0170-4214},
    review={MR2102321},
}



\bib{MR1026740}{article}{
    author={Boisseau, B.},
    author={van Leeuwen, W. A.},
     title={Relativistic Boltzmann theory in $D+1$ spacetime dimensions},
   journal={Ann. Physics},
    volume={195},
      date={1989},
    number={2},
     pages={376\ndash 419},
      issn={0003-4916},
    review={MR1026740 (90j:82019)},
}




\bib{MR0213137}{article}{
    author={Bichteler, Klaus},
     title={On the Cauchy problem of the relativistic Boltzmann equation},
   journal={Comm. Math. Phys.},
    volume={4},
      date={1967},
     pages={352\ndash 364},
      issn={0010-3616},
    review={MR0213137 (35 \#4002)},
}


\bib{MR2098116}{article}{
    author={Calogero, Simone},
     title={The Newtonian limit of the relativistic Boltzmann equation},
   journal={J. Math. Phys.},
    volume={45},
      date={2004},
    number={11},
     pages={4042\ndash 4052},
      issn={0022-2488},
    review={MR2098116},
}



\bib{MR700073}{article}{
   author={Cercignani, Carlo},
   title={Speed of propagation of infinitesimal disturbances in a
   relativistic gas},
   journal={Phys. Rev. Lett.},
   volume={50},
   date={1983},
   number={15},
   pages={1122--1124},
   issn={0031-9007},
   review={\MR{700073 (84f:82027)}},
   doi={10.1103/PhysRevLett.50.1122},
}		



		
\bib{MR1770447}{article}{
   author={Cercignani, Carlo},
   title={Propagation phenomena in classical and relativistic rarefied
   gases},
   booktitle={Proceedings of the Fifth International Workshop on
   Mathematical Aspects of Fluid and Plasma Dynamics (Maui, HI, 1998)},
   journal={Transport Theory Statist. Phys.},
   volume={29},
   date={2000},
   number={3-5},
   pages={607--614},
   issn={0041-1450},
   review={\MR{1770447 (2001f:82058)}},
   doi={10.1080/00411450008205896},
}




\bib{MR1307620}{book}{
   author={Cercignani, Carlo},
   author={Illner, Reinhard},
   author={Pulvirenti, Mario},
   title={The mathematical theory of dilute gases},
   series={Applied Mathematical Sciences},
   volume={106},
   publisher={Springer-Verlag},
   place={New York},
   date={1994},
   pages={viii+347},
   isbn={0-387-94294-7},
   review={\MR{1307620 (96g:82046)}},
}







\bib{MR1706733}{article}{
   author={Cercignani, C.},
   author={Kremer, G. M.},
   title={On relativistic collisional invariants},
   journal={J. Statist. Phys.},
   volume={96},
   date={1999},
   number={1-2},
   pages={439--445},
   issn={0022-4715},
   review={\MR{1706733 (2000h:82037)}},
   doi={10.1023/A:1004545104959},
}


\bib{MR1745845}{article}{
   author={Cercignani, C.},
   author={Kremer, G. M.},
   title={Trend to equilibrium of a degenerate relativistic gas},
   journal={J. Statist. Phys.},
   volume={98},
   date={2000},
   number={1-2},
   pages={441--456},
   issn={0022-4715},
   review={\MR{1745845 (2000m:82048)}},
   doi={10.1023/A:1018695426728},
}
		
		




\bib{MR1810262}{article}{
   author={Cercignani, C.},
   author={Kremer, G. M.},
   title={Moment closure of the relativistic Anderson and Witting model
   equation},
   journal={Phys. A},
   volume={290},
   date={2001},
   number={1-2},
   pages={192--202},
   issn={0378-4371},
   review={\MR{1810262 (2001j:82088)}},
   doi={10.1016/S0378-4371(00)00403-9},
}
		

\bib{MR1857127}{article}{
   author={Cercignani, C.},
   author={Kremer, G. M.},
   title={Dispersion and absorption of plane harmonic waves in a
   relativistic gas},
   journal={Contin. Mech. Thermodyn.},
   volume={13},
   date={2001},
   number={3},
   pages={171--182},
   issn={0935-1175},
   review={\MR{1857127 (2002g:76121)}},
}



		





\bib{MR1898707}{book}{
    author={Cercignani, Carlo},
    author={Kremer, Gilberto Medeiros},
     title={The relativistic Boltzmann equation: theory and applications},
    series={Progress in Mathematical Physics},
    volume={22},
 publisher={Birkh\"auser Verlag},
     place={Basel},
      date={2002},
     pages={x+384},
      isbn={3-7643-6693-1},
    review={MR1898707 (2003f:82078)},
}

\bib{MR767047}{article}{
   author={Cercignani, Carlo},
   author={Majorana, Armando},
   title={Propagation of infinitesimal disturbances in a gas according to a
   relativistic kinetic model},
   language={English, with Italian summary},
   journal={Meccanica},
   volume={19},
   date={1984},
   number={3},
   pages={175--181},
   issn={0025-6455},
   review={\MR{767047 (85k:76045)}},
   doi={10.1007/BF01743729},
}
		




\bib{MR635279}{book}{
    author={de Groot, S. R.},
    author={van Leeuwen, W. A.},
    author={van Weert, Ch. G.},
     title={Relativistic kinetic theory},
 publisher={North-Holland Publishing Co.},
     place={Amsterdam},
      date={1980},
     pages={xvii+417},
      isbn={0-444-85453-3},
    review={MR635279 (83a:82024)},
}





\bib{MR0471665}{article}{
    author={Dijkstra, J. J.},
    author={van Leeuwen, W. A.},
     title={Mathematical aspects of relativistic kinetic theory},
   journal={Phys. A},
    volume={90},
      date={1978},
    number={3--4},
     pages={450\ndash 486},
    review={MR0471665 (57 \#11390)},
}



\bib{MR1014927}{article}{
    author={DiPerna, R. J.},
    author={Lions, P.-L.},
     title={On the Cauchy problem for Boltzmann equations: global existence
            and weak stability},
   journal={Ann. of Math. (2)},
    volume={130},
      date={1989},
    number={2},
     pages={321\ndash 366},
      issn={0003-486X},
    review={MR1014927 (90k:82045)},
}






\bib{DEnotMSI}{article}{
   author={Dudy{\'n}ski, Marek},
   author={Ekiel-Je{\.z}ewska, Maria L.},
   title={The relativistic Boltzmann equation - mathematical and physical aspects},
   journal={J. Tech. Phys.},
   volume={48},
   date={2007},
   pages={39--47},
}



\bib{MR1031410}{article}{
    author={Dudy{\'n}ski, Marek},
     title={On the linearized relativistic Boltzmann equation. II. Existence
            of hydrodynamics},
   journal={J. Statist. Phys.},
    volume={57},
      date={1989},
    number={1-2},
     pages={199\ndash 245},
      issn={0022-4715},
    review={MR1031410 (91b:82043)},
}


\bib{MR933458}{article}{
    author={Dudy{\'n}ski, Marek},
    author={Ekiel-Je{\.z}ewska, Maria L.},
     title={On the linearized relativistic Boltzmann equation. I. Existence
            of solutions},
   journal={Comm. Math. Phys.},
    volume={115},
      date={1988},
    number={4},
     pages={607\ndash 629},
      issn={0010-3616},
    review={MR933458 (89h:82017)},
}


\bib{MR1151987}{article}{
    author={Dudy{\'n}ski, Marek},
    author={Ekiel-Je{\.z}ewska, Maria L.},
     title={Global existence proof for relativistic Boltzmann equation},
   journal={J. Statist. Phys.},
    volume={66},
      date={1992},
    number={3-4},
     pages={991\ndash 1001},
      issn={0022-4715},
    review={MR1151987 (93b:82064)},
}

\bib{MR841735}{article}{
    author={Dudy{\'n}ski, Marek},
    author={Ekiel-Je{\.z}ewska, Maria L.},
     title={Errata: ``Causality of the linearized relativistic Boltzmann
            equation''},
   journal={Investigaci\'on Oper.},
    volume={6},
      date={1985},
    number={1},
     pages={2228},
      issn={0257-4306},
    review={MR841735 (87e:82043b)},
}

\bib{MR818441}{article}{
    author={Dudy{\'n}ski, Marek},
    author={Ekiel-Je{\.z}ewska, Maria L.},
     title={Causality of the linearized relativistic Boltzmann equation},
   journal={Phys. Rev. Lett.},
    volume={55},
      date={1985},
    number={26},
     pages={2831\ndash 2834},
      issn={0031-9007},
    review={MR818441 (87e:82043a)},
}



\bib{MR2459827}{article}{
   author={Ha, Seung-Yeal},
   author={Kim, Yong Duck},
   author={Lee, Ho},
   author={Noh, Se Eun},
   title={Asymptotic completeness for relativistic kinetic equations with
   short-range interaction forces},
   journal={Methods Appl. Anal.},
   volume={14},
   date={2007},
   number={3},
   pages={251--262},
   issn={1073-2772},
   review={\MR{2459827 (2010b:35464)}},
}

\bib{MR2543323}{article}{
   author={Ha, Seung-Yeal},
   author={Lee, Ho},
   author={Yang, Xiongfeng},
   author={Yun, Seok-Bae},
   title={Uniform $L^2$-stability estimates for the relativistic
   Boltzmann equation},
   journal={J. Hyperbolic Differ. Equ.},
   volume={6},
   date={2009},
   number={2},
   pages={295--312},
   issn={0219-8916},
   review={\MR{2543323}},
   doi={10.1142/S0219891609001848},
}

\bib{MR2249574}{article}{
   author={Hsiao, Ling},
   author={Yu, Hongjun},
   title={Asymptotic stability of the relativistic Maxwellian},
   journal={Math. Methods Appl. Sci.},
   volume={29},
   date={2006},
   number={13},
   pages={1481--1499},
   issn={0170-4214},
   review={\MR{2249574 (2008f:35393)}},
   doi={10.1002/mma.736},
}

\bib{MR2289548}{article}{
   author={Hsiao, Ling},
   author={Yu, Hongjun},
   title={Global classical solutions to the initial value problem for the
   relativistic Landau equation},
   journal={J. Differential Equations},
   volume={228},
   date={2006},
   number={2},
   pages={641--660},
   issn={0022-0396},
   review={\MR{2289548 (2007k:35504)}},
   doi={10.1016/j.jde.2005.10.022},
}
		


\bib{MR1379589}{book}{
   author={Glassey, Robert T.},
   title={The Cauchy problem in kinetic theory},
   publisher={Society for Industrial and Applied Mathematics (SIAM)},
   place={Philadelphia, PA},
   date={1996},
   pages={xii+241},
   isbn={0-89871-367-6},
   review={\MR{1379589 (97i:82070)}},
}

\bib{MR2217287}{article}{
   author={Glassey, Robert T.},
   title={Global solutions to the Cauchy problem for the relativistic
   Boltzmann equation with near-vacuum data},
   journal={Comm. Math. Phys.},
   volume={264},
   date={2006},
   number={3},
   pages={705--724},
   issn={0010-3616},
   review={\MR{2217287 (2007a:82062)}},
}



\bib{MR1105532}{article}{
    author={Glassey, Robert T.},
    author={Strauss, Walter A.},
     title={On the derivatives of the collision map of relativistic
            particles},
   journal={Transport Theory Statist. Phys.},
    volume={20},
      date={1991},
    number={1},
     pages={55\ndash 68},
      issn={0041-1450},
    review={MR1105532 (92f:81222)},
}

\bib{MR1211782}{article}{
    author={Glassey, Robert T.},
    author={Strauss, Walter A.},
     title={Asymptotic stability of the relativistic Maxwellian},
   journal={Publ. Res. Inst. Math. Sci.},
    volume={29},
      date={1993},
    number={2},
     pages={301\ndash 347},
      issn={0034-5318},
    review={MR1211782 (94c:82063)},
}

\bib{MR1321370}{article}{
    author={Glassey, Robert T.},
    author={Strauss, Walter A.},
     title={Asymptotic stability of the relativistic Maxwellian via fourteen
            moments},
   journal={Transport Theory Statist. Phys.},
    volume={24},
      date={1995},
    number={4-5},
     pages={657\ndash 678},
      issn={0041-1450},
    review={MR1321370 (96c:82054)},
}







\bib{MR2000470}{article}{
   author={Guo, Yan},
   title={The Vlasov-Maxwell-Boltzmann system near Maxwellians},
   journal={Invent. Math.},
   volume={153},
   date={2003},
   number={3},
   pages={593--630},
   issn={0020-9910},
   review={\MR{2000470 (2004m:82123)}},
}



\bib{gsRVMB}{article}{
   author={Guo, Yan},
   author={Strain, Robert M.},
   title={Momentum Regularity and Stability of the Relativistic Vlasov-Maxwell-Boltzmann System},
   volume={in preparation},
   date={2010},
}







\bib{MR1714446}{article}{
    author={Jiang, Zhenglu},
     title={On the Cauchy problem for the relativistic Boltzmann equation in
            a periodic box: global existence},
   journal={Transport Theory Statist. Phys.},
    volume={28},
      date={1999},
    number={6},
     pages={617\ndash 628},
      issn={0041-1450},
    review={MR1714446 (2000k:82076)},
}

\bib{MR1676150}{article}{
    author={Jiang, Zhenglu},
     title={On the relativistic Boltzmann equation},
   journal={Acta Math. Sci. (English Ed.)},
    volume={18},
      date={1998},
    number={3},
     pages={348\ndash 360},
      issn={0252-9602},
    review={MR1676150 (2001a:82062)},
}


\bib{MR0004796}{article}{
    author={Lichnerowicz, Andr{\'e}},
    author={Marrot, Raymond},
     title={Propri\'et\'es statistiques des ensembles de particules en
            relativit\'e restreinte},
  language={French},
   journal={C. R. Acad. Sci. Paris},
    volume={210},
      date={1940},
     pages={759\ndash 761},
      issn={0249-6305},
    review={MR0004796 (3,63b)},
}

\bib{MR1284432}{article}{
   author={Lions, P.-L.},
   title={Compactness in Boltzmann's equation via Fourier integral operators
   and applications. I, II, III},
   journal={J. Math. Kyoto Univ.},
   volume={34},
   date={1994},
   number={2,3},
   pages={391--427, 429--461,539--584},
   issn={0023-608X},
   review={\MR{1284432 (96f:35133)}},
}




\bib{ssHilbert}{article}{
   author={Speck, Jared},
      author={Strain, Robert M.},
   title={Hilbert Expansion from the Boltzmann equation to relativistic Fluids},
   journal={Comm. Math. Phys.},
   volume={accepted},
   date={2010},
   pages={49pp},
}



\bib{MR0088362}{book}{
    author={Synge, J. L.},
     title={The relativistic gas},
 publisher={North-Holland Publishing Company, Amsterdam},
      date={1957},
     pages={xi+108},
    review={MR0088362 (19,508b)},
}

\bib{stewart}{book}{
    author={Stewart, J. M.},
     title={Non-equilibrium relativistic kinetic theory, volume 10 of Lectures notes in physics},
 publisher={Springer-Verlag, Berlin},
      date={1971},
}

\bib{MR2100057}{article}{
   author={Strain, Robert M.},
   author={Guo, Yan},
   title={Stability of the relativistic Maxwellian in a collisional plasma},
   journal={Comm. Math. Phys.},
   volume={251},
   date={2004},
   number={2},
   pages={263--320},
   review={\MR{2100057 (2005m:82155)}},
}



\bib{strainPHD}{book}{
    author={Strain, Robert M.},
     title={An Energy Method in Collisional Kinetic Theory,},
 publisher={Ph.D. dissertation, Division of Applied Mathematics, Brown University},
      date={May 2005},
}





\bib{strainNEWT}{article}{
    author={Strain, Robert M.},
     title={Global Newtonian Limit For The  Relativistic
Boltzmann Equation Near Vacuum},
   journal={SIAM J. MATH. ANAL.},
   volume={42},
   date={2010},
   number={4},
   pages={1568Ð1601},
   doi={110.1137/090762695},
}


\bib{strainSOFT}{article}{
    author={Strain, Robert M.},
     title={Asymptotic Stability of the Relativistic {B}oltzmann Equation for the Soft-Potentials},
   journal={Comm. Math. Phys.},
   volume={300},
   date={2010},
   number={2},
   pages={529--597},
   eprint={arXiv:1003.4893v1}
      doi={10.1007/s00220-010-1129-1},
}



\bib{weinbergBK}{book}{
    author={Weinberg, Stephen},
     title={Gravitation and Cosmology: Principles and Applications of the General Theory of Relativity},
 publisher={Wiley, New York},
      date={1972},
}

\bib{MR1480243}{article}{
    author={Wennberg, Bernt},
     title={The geometry of binary collisions and generalized Radon
            transforms},
   journal={Arch. Rational Mech. Anal.},
    volume={139},
      date={1997},
    number={3},
     pages={291\ndash 302},
      issn={0945-8396},
    review={MR1480243 (98k:82166)},
}

\bib{MR2593052}{article}{
   author={Yang, Tong},
   author={Yu, Hongjun},
   title={Hypocoercivity of the relativistic Boltzmann and Landau equations
   in the whole space},
   journal={J. Differential Equations},
   volume={248},
   date={2010},
   number={6},
   pages={1518--1560},
   issn={0022-0396},
   review={\MR{2593052}},
   doi={10.1016/j.jde.2009.11.027},
}

\bib{MR2514726}{article}{
   author={Yu, Hongjun},
   title={Smoothing effects for classical solutions of the relativistic
   Landau-Maxwell system},
   journal={J. Differential Equations},
   volume={246},
   date={2009},
   number={10},
   pages={3776--3817},
   issn={0022-0396},
   review={\MR{2514726 (2010e:35282)}},
   doi={10.1016/j.jde.2009.02.021},
}
		


\end{biblist}
\end{bibdiv}
\end{document}